\documentclass{article}
\usepackage{latexsym,amsthm,amssymb,amscd,amsmath}
\newcounter{alphthm}
\newtheorem{thm}{Theorem}

\newtheorem{prop}{Proposition}
\newtheorem{cor}{Corollary}
\newtheorem{lem}{Lemma}

\newcommand{\be}{\begin{equation}}
\newcommand{\ee}{\end{equation}}
\newcommand{\ben}{\begin{enumerate}}
\newcommand{\een}{\end{enumerate}}
\newcommand{\pa}{{\partial}}

\def\beq{\begin{equation}}
\def\eeq{\end{equation}}

\title{Doubly Warped Product Finsler Manifolds with Some Non-Riemannian Curvature Properties}
\author{E. Peyghan, A. Tayebi and B. Najafi}
\begin{document}

\maketitle
\begin{abstract}
In this paper, we consider doubly warped product (DWP) Finsler manifolds with some non-Riemannian curvature properties. First, we study Berwald and isotropic mean Berwald DWP-Finsler manifolds. Then we prove that every proper  Douglas DWP-Finsler manifold  is Riemannian. We show that a proper DWP-manifold is Landsbergian if and only it is Berwaldian. In continue, we prove that every relatively isotropic Landsberg DWP-manifold is a Landsberg manifold. We show that a relatively isotropic mean Landsberg warped product manifold is a weakly Landsberg manifold. Finally, we  show that there is no any locally dually flat proper DWP-Finsler manifold. 
\end{abstract}

\section{Introduction}

The study of relativity theory demands a wider class of manifold and the idea of doubly warped products was introduced and studied by many authors. The recently studies show that  the notion of doubly warped product manifolds has an important role in Riemannian geometry and its application \cite{A}\cite{BEP}\cite{BP}\cite{G}\cite{Marian}\cite{Marian1} \cite{U}. For example, Beem-Powell study this product for Lorentzian manifolds \cite{BP}.  Then in  \cite{A},  Allison considered  global hyperbolicity of doubly warped product and null pseudo convexity of Lorentizian doubly warped product.

On the other hand, Finsler geometry consists of classical and generalized Finsler geometries. It is a subject studying manifolds
whose tangent spaces carry a norm varying smoothly with the base point. Indeed,  Finsler geometry is just
Riemannian geometry without the quadratic restriction. Thus it is  natural to extending the construction of  warped product manifolds for Finsler geometry. In the first step,  Asanov  gave the generalization of the Schwarzschild metric in the Finslerian setting and  obtain some models of relativity theory described through the warped product of Finsler metrics \cite{As}\cite{As1}. In  \cite{Koz},  Kozma-Peter-Varga define their warped product for Finsler metrics and  conclude that completeness of doubly warped product can be related to completeness of  its components.

Let $(M_1,F_1)$ and $(M_2,F_2)$ be two Finsler manifolds  and $f_1: M_1\rightarrow \mathbb{R}^+$ and $f_2: M_2\rightarrow \mathbb{R}^+$ are two smooth functions. Let $\pi_1: M_1\times M_2\rightarrow M_1$ and $\pi_2: M_1\times M_2\rightarrow M_2$ be the natural projection maps. The product manifold $M_1\times M_2$ endowed with the metric $F: TM^\circ_1\times TM^\circ_2\rightarrow\mathbb{R}$ is considered,
\begin{equation}
F(v_1, v_2)=\sqrt{f_2^2(\pi_2(v_2))F_1^2(v_1)+f_1^2(\pi_1(v_1))F_2^2(v_2)},\label{Finsler}
\end{equation}
where $TM^\circ_1=TM_1-\{0\}$ and $TM^\circ_2=TM_2-\{0\}$. The metric defined above is a Finsler metric. The product manifold $M_1\times M_2$ with the metric $F(v)=F(v_1,v_2)$ for $(v_1,v_2)\in TM_1^\circ\times TM_2^\circ$ defined above will be called the doubly warped product of the manifolds $M_1$ and $M_2$ and $f_1$ and $f_2$ will be called the warping functions. We denote this doubly warped by ${}_{f_2}M_1\times_{f_1}M_2$. If either $f_1=1$ or $f_2=1$, then ${}_{f_2}M_1\times_{f_1}M_2$ becomes a warped product of Finsler manifolds $M_1$ and $M_2$. If $f_1=f_2=1$, then we have a product manifold. If neither $f_1$ nor $f_2$ is constant, then we have a proper DWP-manifold.

Let $(M, F)$ be a Finsler manifold.  The second and third order derivatives of ${1\over 2} F_x^2$ at $y\in T_xM_{0}$ are the symmetric trilinear forms ${\bf g}_y$ and ${\bf C}_y$ on $T_xM$, which called the fundamental tensor and   Cartan torsion, respectively. The rate of change of ${\bf C}_y$ along geodesics is the
Landsberg curvature  ${\bf L}_y$ on $T_xM$ \cite{Ba}\cite{BCS}.   $F$ is said to be  relatively isotropic Landsberg metric if it satisfies  ${\bf L}+cF{\bf C}=0$, where $c=c(x)$ is a scalar function on $M$. Set ${\bf I}_y:= \sum_{i=1}^n {\bf C}_y(e_i, e_i, \cdot )$ and ${\bf J}_y:= \sum_{i=1}^n {\bf L}_y(e_i, e_i, \cdot )$, where $\{e_i\}$ is an orthonormal basis for $(T_xM, {\bf g}_y)$. ${\bf I}_y$ and  ${\bf J}_y$ is called the mean Cartan torsion and  mean Landsberg curvature, respectively.  $F$  is said to be relatively isotropic mean Landsberg metric if ${\bf J}+cF {\bf I}=0$, where $c=c(x)$ is a scalar function on $M$ \cite{ChSh}.

The geodesic curves of a Finsler metric $F$ on a smooth manifold $M$, are determined  by the system of second order differential equations $ \ddot c^i+2G^i(\dot c)=0$, where the local functions $G^i=G^i(x, y)$ are called the  spray coefficients.  $F$ is  called a Berwald metric, if  $G^i$  are quadratic in $y\in T_xM$  for any $x\in M$,  and called a Douglas metric if  $G^i=\frac{1}{2}\Gamma^i_{jk}(x)y^jy^k + P(x,y)y^i$ \cite{BM2}\cite{NST1}. Taking a trace of  Berwald curvature yields mean Berwald curvature ${\bf E}$.  Then $F$ is said to be isotropic mean Berwald metric if ${\bf E}=\frac{n+1}{2}cF^{-1}{\bf h}$, where ${\bf h}=h_{ij} dx^i\otimes dx^j$ is the angular metric and $c=c(x)$ is a scalar function on $M$ \cite{NST2}.

This paper is arranged as follows: In section 2, we give some  basic concepts of Finsler manifolds. In sections 3 and 4, we study doubly warped product Finsler metrics or briefly DWP-Finsler metrics with vanishing Berwald curvature and isotropic mean Berwald curvature, respectively. In section 5, we prove that every proper Douglas  DWP-Finsler manifold  is Riemannian. In section 6, we show that a proper DWP-Finsler manifold is a Landsberg manifold if and only if it is a Berwald manifold. Then we prove that every relatively isotropic Landsberg DWP-Finsler manifold is a Landsberg manifold. In section 7, we prove that a relatively isotropic mean Landsberg warped product manifold is a weakly Landsberg manifold. Finally in section 8, we show that there is no any locally dually flat proper DWP-Finsler manifold.

\section{Preliminaries}
Let $M$ be an $n$-dimensional $ C^\infty$ manifold. Denote by $T_x M $ the tangent space at $x \in M$,  by $TM=\cup _{x \in M} T_x M $ the tangent bundle of $M$, and by $TM_{0} = TM \setminus \{ 0 \}$ the slit tangent bundle on $M$. A  Finsler metric on $M$ is a function $ F:TM \rightarrow [0,\infty)$ which has the following properties:\\
(i) $F$ is $C^\infty$ on $TM_{0}$;\\
(ii) $F$ is positively 1-homogeneous on the fibers of tangent bundle $TM$;\\
(iii) for each $y\in T_xM$, the following quadratic form ${\bf g}_y$ on
$T_xM$  is positive definite,
\[
{\bf g}_{y}(u,v):={1 \over 2} \frac{\partial^2}{\partial s \partial t}\left[  F^2 (y+su+tv)\right]|_{s,t=0}, \ \
u,v\in T_xM.
\]
Let  $x\in M$ and $F_x:=F|_{T_xM}$.  To measure the non-Euclidean feature of $F_x$, define ${\bf C}_y:T_xM\otimes T_xM\otimes T_xM\rightarrow \mathbb{R}$ by
\[
{\bf C}_{y}(u,v,w):={1 \over 2} \frac{d}{dt}\left[{\bf g}_{y+tw}(u,v)
\right]|_{t=0}, \ \ u,v,w\in T_xM.
\]
The family ${\bf C}:=\{{\bf C}_y\}_{y\in TM_0}$  is called the Cartan torsion. It is well known that ${\bf{C}}=0$ if and only if $F$ is Riemannian. For $y\in T_x M_0$, define  mean Cartan torsion ${\bf I}_y$ by ${\bf I}_y(u):=I_i(y)u^i$, where $I_i:=g^{jk}C_{ijk}$, $C_{ijk}=\frac{1}{2}\frac{\pa g_{ij}}{\pa y^k}$ and $u=u^i\frac{\partial}{\partial x^i}|_x$. By Deicke's  Theorem, $F$ is Riemannian  if and only if ${\bf I}_y=0$ \cite{BCS}\cite{ShDiff}.

\bigskip

Given a Finsler manifold $(M,F)$, then a global vector field ${\bf G}$ is induced by $F$ on $TM_0$, which in a standard coordinate $(x^i,y^i)$ for $TM_0$ is given by ${\bf G}=y^i {{\partial} \over {\partial x^i}}-2G^i(x,y){{\partial} \over {\partial y^i}}$, where \[
G^i:=\frac{1}{4}g^{il}(y)\Big\{\frac{\partial^2 F^2}{\partial x^k
\partial y^l}y^k-\frac{\partial F^2}{\partial x^l}\Big\},\ \ \ y\in T_xM.
\]
{\bf G} is called the  spray associated  to $(M,F)$.  In local coordinates, a curve $c(t)$ is a geodesic if and only if its coordinates $(c^i(t))$ satisfy $ \ddot c^i+2G^i(\dot c)=0$.
\bigskip

For a tangent vector $y \in T_xM_0$, define ${\bf B}_y:T_xM\otimes T_xM \otimes T_xM\rightarrow T_xM$ and ${\bf E}_y:T_xM \otimes T_xM\rightarrow \mathbb{R}$ by ${\bf B}_y(u, v, w):=B^i_{\ jkl}(y)u^jv^kw^l{{\partial } \over {\partial x^i}}|_x$ and ${\bf E}_y(u,v):=E_{jk}(y)u^jv^k$
where
\[
B^i_{\ jkl}:={{\partial^3 G^i} \over {\partial y^j \partial y^k \partial y^l}},\ \ \ E_{jk}:=\frac{1}{2}B^m_{\ jkm}.
\]
$\bf B$ and $\bf E$ are called the Berwald curvature and mean Berwald curvature, respectively.  Then $F$ is called a Berwald metric and weakly Berwald metric if ${\bf{B}}=0$ and ${\bf{E}}=0$, respectively  \cite{ShDiff}\cite{TP}. It is  proved that on  a Berwald space,  the parallel translation along any geodesic preserves the Minkowski functionals \cite{Ich}. Thus Berwald spaces can be viewed as Finsler spaces modeled on a single Minkowski space.

A Finsler metric $F$ is said to be isotropic mean Berwald metric if its  mean Berwald curvature is in the following form
\begin{eqnarray}\label{IBCurve}
E_{ij}=\frac{1}{2}(n+1)cF^{-1}h_{ij},
\end{eqnarray}
where $h_{ij}=g_{ij}-F^{-2}y_iy_j$ is the angular metric and $c=c(x)$ is a  scalar function  on $M$ \cite{ChSh}.

\bigskip

Define ${\bf D}_y:T_xM\otimes T_xM \otimes T_xM\rightarrow T_xM$  by
${\bf D}_y(u,v,w):=D^i_{\ jkl}(y)u^iv^jw^k\frac{\partial}{\partial x^i}|_{x}$ where
\[
D^i_{\ jkl}:=B^i_{\ jkl}-{2\over n+1}\{E_{jk}\delta^i_l+E_{jl}\delta^i_k+E_{kl}\delta^i_j+\frac{E_{jk}}{\partial y^l}y^i\}.
\]
We call ${\bf D}:=\{{\bf D}_y\}_{y\in TM_{0}}$ the Douglas curvature. A Finsler metric with ${\bf D}=0$ is called a Douglas metric. The notion of Douglas metrics was proposed by B$\acute{a}$cs$\acute{o}$-Matsumoto as a generalization  of Berwald metrics \cite{BM2}.

\bigskip

There is another extension of Berwald curvature. For a tangent vector $y \in T_xM_0$, define  ${\bf L}_y:T_xM\otimes T_xM\otimes T_xM\rightarrow \mathbb{R}$  by ${\bf L}_y(u,v,w):=L_{ijk}(y)u^iv^jw^k$, where
\[
L_{ijk}:=\frac{-1}{2}y_l B^l_{\ ijk}.
\]
The family ${\bf L}:=\{{\bf L}_y\}_{y\in TM_{0}}$  is called the Landsberg curvature. A Finsler metric is called a Landsberg metric  if ${\bf{L}}=0$. The quantity ${\bf L}/{\bf C}$ is regarded as the relative rate of change of ${\bf C}$ along geodesics.  A Finsler metric $F$  is said to be  relatively isotropic Landsberg metric if it satisfies
\[
{\bf L}=cF \bf C,
\]
where $c=c(x)$ is a scalar function on $M$ \cite{ChSh}.

\smallskip

Taking a trace of  Landsberg curvature yields mean   Landsberg curvature ${\bf J}_y:T_xM\rightarrow \mathbb{R}$, which defined by  ${\bf J}_y(u): = J_i (y)u^i$, where
\[
J_i: = g^{jk}L_{ijk}.
\]
A Finsler metric is called a weakly  Landsberg metric  if ${\bf{J}}=0$. The quantity ${\bf J}/{\bf I}$ is regarded as the relative rate of change of ${\bf I}$ along geodesics.  A Finsler metric $F$  is said to be relatively isotropic mean Landsberg metric if
\[
{\bf J}=cF \bf I,
\]
for some scalar function $c=c(x)$ on $M$ \cite{ChSh}. It is obvious that,  every  relatively isotropic Landsberg metric is a  relatively isotropic mean  Landsberg metric.

\bigskip

 A Finsler metric $F=F(x,y)$ on a manifold $M$ is said to be locally dually flat if at any point there is a standard coordinate system
 $(x^i,y^i)$ in $TM$  satisfies
\be
(F^2)_{x^ky^l}y^k=2(F^2)_{x^l}.\label{LDF}
\ee
In this case, the coordinate $(x^i)$ is called an adapted local coordinate system \cite{am}\cite{amna}. It is easy to see that, every locally Minkowskian metric satisfies (\ref{LDF}), hence is locally dually flat \cite{shenInf}. But the converse is not true, generally.

\section{Berwaldian DWP-Finsler Manifolds}
In this  section, we study  DWP-Finsler manifolds with  vanishing Berwald curvature.
\begin{lem}\label{Lemma1}
Every proper  DWP-Finsler manifold $({}_{f_2}M_1\times{}_{f_1}M_2,F)$ with vanishing Berwald curvature is a Riemannian manifold.
\end{lem}
\begin{proof}
The Berwald curvature of a DWP-Finsler manifold $({}_{f_2}M_1\times{}_{f_1}M_2,F)$ are given as follows:
\begin{eqnarray}
&&\textbf{B}^k_{ijl}=B^k_{ijl}-\frac{1}{4f_2^2}\frac{\partial^3g^{kh}}{\partial y^i\partial y^j\partial y^l}\frac{\partial f_1^2}{\partial x^h}F_2^2,\label{feri}
\\
&&\textbf{B}^k_{i\beta l}=-\frac{1}{4f_2^2}\frac{\partial^2g^{kh}}{\partial y^l\partial y^i}\frac{\partial f_1^2}{\partial x^h}\frac{\partial F_2^2}{\partial v^\beta},
\\
&&\textbf{B}^k_{\alpha\beta l}=-\frac{1}{f_2^2}\frac{\partial f_1^2}{\partial x^h} \frac{\partial g^{kh}}{\partial y^l}g_{\alpha\beta},
\\
&&\textbf{B}^k_{\alpha\beta\lambda}=-\frac{1}{f_2^2}\frac{\partial f_1^2}{\partial x^h}g^{kh}C_{\alpha\beta\lambda},\label{1b}
\\
&&\textbf{B}^\gamma_{\alpha\beta\lambda}=B^\gamma_{\alpha\beta\lambda}-\frac{1}{4f_1^2}
\frac{\partial^3g^{\gamma\nu}}{\partial v^\beta\partial v^\alpha\partial v^\lambda}\frac{\partial f_2^2}{\partial u^\nu}F_1^2,\label{4}\\
&&\textbf{B}^\gamma_{i\beta\lambda}=-\frac{1}{4f_1^2}\frac{\partial^2g^{\alpha\gamma}}{\partial v^\beta\partial v^\lambda}\frac{\partial f_2^2}{\partial u^\alpha}\frac{\partial F_1^2}{\partial y^i},\\
&&\textbf{B}^\gamma_{ij\lambda}=-\frac{1}{2f_1^2}g_{ij}\frac{\partial g^{\alpha\gamma}}{\partial v^\lambda}\frac{\partial f_2^2}{\partial u^\alpha},\label{3}\\
&&\textbf{B}^\gamma_{ijk}=-\frac{1}{f_1^2}C_{ijk}g^{\alpha\gamma}\frac{\partial f_2^2}{\partial u^\alpha}\label{2a},
\end{eqnarray}
If $({}_{f_2}M_1\times{}_{f_1}M_2,F)$ is Berwaldian, then we have $\textbf{B}^d_{abc}=0$. By (\ref{1b}), we  get
\be
C_{\alpha\beta\lambda}g^{kh}\frac{\partial f_1^2}{\partial x^h}=0.\label{B1}
\ee
Multiplying (\ref{B1}) with $g_{kr}$ implies that
\be
C_{\alpha\beta\lambda}\frac{\partial f_1^2}{\partial x^r}=0.\label{B2}
\ee
By  (\ref{B2}),  if $f_1$ is not constant then we get $C_{\alpha\beta\lambda}=0$, i.e., $(M_2, F_2)$ is Riemannian. In the similar way, from (\ref{2a}) we conclude that if $f_2$ is non constant then $(M_1, F_1)$ is Riemannian.
\end{proof}
\begin{thm}\label{Theorem1}
Let $({}_{f_2}M_1\times{}_{f_1}M_2,F)$ be a DWP-Finsler manifold.\\
(i) If $f_1$ is constant and $f_2$ is not constant. Then $({}_{f_2}M_1\times{}_{f_1}M_2,F)$ is a Berwald manifold if and only if $M_1$ is Riemannian, $M_2$ is a Berwald manifold and
\[
\frac{\partial g^{\alpha\gamma}}{\partial v^\lambda}\frac{\partial f_2^2}{\partial u^\alpha}=0.
\]
(ii) If $f_2$ is constant and $f_1$ is not constant. Then $({}_{f_2}M_1\times{}_{f_1}M_2,F)$ is a Berwald manifold if and only if $M_2$ is Riemannian, $M_1$ is Berwaldian and
\[
\frac{\partial g^{ij}}{\partial y^k}\frac{\partial f_1^2}{\partial x^i}=0.
\]
\end{thm}
\begin{proof} Let $({}_{f_2}M_1\times{}_{f_1}M_2,F)$ be a Berwaldian manifold and $f_1$ is constant on $M_1$. Then  from (\ref{2a}) we get $C_{ijk}=0$, i.e., $(M_1, F_1)$ is Riemannian. Also, (\ref{3}) gives us
\[
\frac{\partial g^{\alpha\gamma}}{\partial v^\lambda}\frac{\partial f_2^2}{\partial u^\alpha}=0.
\]
Differentiating this equation with respect to $v^{\beta}$, implies that
\[
\frac{\partial^2 g^{\alpha\gamma}}{\partial v^\lambda\partial v^\beta}\frac{\partial f_2^2}{\partial u^\alpha}=0,
\]
and consequently
\[
\frac{\partial^3 g^{\alpha\gamma}}{\partial v^\lambda\partial v^\beta\partial v^\mu}\frac{\partial f_2^2}{\partial u^\alpha}=0.
\]
Then (\ref{4}) reduces to  $B^\gamma_{\alpha\beta\lambda}=0$, i.e., $(M_2, F_2)$ is Berwaldian. In the similar way, we can prove the converse of this assertion. The proof of (ii) is similar to that of (i) and we omit it here.
\end{proof}

\smallskip

By using the similar  argument used in Theorem \ref{Theorem1}, we have the following.
\begin{cor}
Let $(M_1\times{}_{f_1}M_2,F)$ be a proper WP-Finsler manifold. Then $(M_1\times{}_{f_1}M_2,F)$ is Berwaldian if and only if $M_2$ is Riemannian, $M_1$ is Berwaldian and
\be
C^{ij}_{\ k}\frac{\partial f_1}{\partial x^i}=0,\label{10}
\ee
where $ C^{ij}_{\ k}=-2\frac{\partial g^{ij}}{\partial y^k}$ is the Cartan tensor.
\end{cor}

\section{Isotropic Mean Berwald DWP-Manifolds}
In this section, we study  DWP-Finsler metrics with isotropic mean Berwald curvature. First, we compute the mean Berwald curvature of a DWP-Finsler
manifold.
\begin{lem}\label{Lemma2}
Let $({}_{f_2}M_1\times{}_{f_1}M_2,F)$ be a DWP-Finsler manifold. Then the mean Berwald curvature of $F$ are given as follows
\begin{eqnarray}
\textbf{E}_{\alpha\beta}\!\!\!\!&=&\!\!\!\!\ E_{\alpha\beta}-\frac{1}{8f_1^2}\frac{\partial^3g^{\gamma\nu}}{\partial v^\beta\partial v^\alpha\partial v^\gamma}\frac{\partial f_2^2}{\partial u^\nu}F_1^2-\frac{1}{4f_2^2}g_{\alpha\beta}\frac{\partial g^{kh}}{\partial y^k}\frac{\partial f_1^2}{\partial x^h},\label{mean1}\\
\textbf{E}_{ij}\!\!\!\!&=&\!\!\!\!\ E_{ij}-\frac{1}{8f_2^2}\frac{\partial^3g^{kh}}{\partial
y^i\partial y^j\partial y^k}\frac{\partial f_1^2}{\partial
x^h}F_2^2-\frac{1}{4f_1^2}g_{ij}\frac{\partial g^{\alpha\gamma}}{\partial v^\gamma}\frac{\partial f_2^2}{\partial u^\alpha},\label{mean2}\\
\textbf{E}_{i\beta
}\!\!\!\!&=&\!\!\!\!-\frac{1}{4f_2^2}\frac{\partial^2g^{kh}}{\partial y^k\partial
y^i}\frac{\partial f_1^2}{\partial x^h}v_\beta-\frac{1}{4f_1^2}\frac{\partial^2g^{\alpha\gamma}}{\partial v^\beta\partial
v^\gamma}\frac{\partial f_2^2}{\partial u^\alpha}y_i,\label{mean3}
\end{eqnarray}
where $E_{ij}$ and $E_{\alpha\beta}$ are mean Berwald curvatures of $(M_1, F_1)$ and $(M_2, F_2)$, respectively.
\end{lem}

\smallskip

\begin{thm}
Let $({}_{f_2}M_1\times{}_{f_1}M_2,F)$ be a proper DWP-Finsler manifold. Then $F$ is a weakly Berwald metric if and only if $F_1$ and $F_2$ are weakly Berwald metric,  and the following holds
\begin{equation}
\frac{\partial g^{kh}}{\partial y^k}\frac{\partial f_1^2}{\partial x^h}=\frac{\partial g^{\gamma\nu}}{\partial v^\gamma}\frac{\partial f_2^2}{\partial u^\nu}=0.\label{equ}
\end{equation}
\end{thm}
\begin{proof}
Let $({}_{f_2}M_1\times_{f_1}M_2, F)$ be a weakly Berwald manifold. Then we have $\textbf{E}_{\alpha\beta}=\textbf{E}_{ij}=\textbf{E}_{i\beta}=0$. Using (\ref{mean3}), we get
\begin{equation}
\frac{1}{f_2^2}\frac{\partial^2g^{kh}}{\partial y^j\partial
y^k}\frac{\partial f_1^2}{\partial x^h}v_\beta=-\frac{1}{f_1^2}\frac{\partial^2g^{\alpha\gamma}}{\partial v^\beta\partial
v^\gamma}\frac{\partial f_2^2}{\partial u^\alpha}y_j.\label{WB1}
\end{equation}
Contracting (\ref{WB1}) with $y^j$ gives us
\be
\frac{1}{f_2^2}\frac{\partial g^{kh}}{\partial y^k}\frac{\partial f_1^2}{\partial x^h}v_\beta=\frac{1}{f_1^2}\frac{\partial^2g^{\nu\gamma}}{\partial v^\beta\partial
v^\gamma}\frac{\partial f_2^2}{\partial u^\nu}F_1^2.\label{WB2}
\ee
Differentiating (\ref{WB2})  with respect $v^\alpha$ implies that
\begin{equation}
\frac{1}{f_2^2}\frac{\partial g^{kh}}{\partial y^k}\frac{\partial f_1^2}{\partial x^h}g_{\alpha\beta}=\frac{1}{f_1^2}\frac{\partial^3g^{\nu\gamma}}{\partial v^\alpha\partial v^\beta\partial
v^\gamma}\frac{\partial f_2^2}{\partial u^\nu}F_1^2.\label{mean4}
\end{equation}
In the similar way, one can obtain
\begin{equation}
\frac{1}{f_1^2}\frac{\partial g^{\gamma\alpha}}{\partial v^\gamma}\frac{\partial f_2^2}{\partial u^\alpha}g_{ij}=\frac{1}{f_2^2}\frac{\partial^3g^{kh}}{\partial y^i\partial y^j\partial
y^k}\frac{\partial f_1^2}{\partial x^h}F_2^2.\label{mean5}
 \end{equation}
Substituting (\ref{mean4}) into (\ref{mean1}) and plugging (\ref{mean5}) into (\ref{mean2}), we have
\begin{equation}
E_{\alpha\beta}=\frac{3}{8f_1^2}\frac{\partial^3g^{\nu\gamma}}{\partial v^\alpha\partial v^\beta\partial
v^\gamma}\frac{\partial f_2^2}{\partial u^\nu}F_1^2,\label{mean6}
\end{equation}
\begin{equation}
E_{ij}=\frac{3}{8f_2^2}\frac{\partial^3g^{kh}}{\partial y^i\partial y^j\partial
y^k}\frac{\partial f_1^2}{\partial x^h}F_2^2.\label{mean7}
\end{equation}
Since $E_{\alpha\beta}$ is a function with respect $(u^\alpha,v^\alpha)$, then by differentiating (\ref{mean6}) with respect $y^h$ we deduce that
\[
\frac{\partial^3g^{\nu\gamma}}{\partial v^\alpha\partial v^\beta\partial
v^\gamma}\frac{\partial f_2^2}{\partial u^\nu}y_h=0,
\]
and consequently
\begin{equation}
\frac{\partial^3g^{\nu\gamma}}{\partial v^\alpha\partial v^\beta\partial
v^\gamma}\frac{\partial f_2^2}{\partial u^\nu}=0.\label{mean8}
\end{equation}
Putting (\ref{mean8})  into (\ref{mean6}) gives us $E_{\alpha\beta}=0$. The similar argument implies that $E_{ij}=0$. Further, from (\ref{mean8}) and (\ref{mean4}) we derive that
\[
\frac{\partial g^{kh}}{\partial y^k}\frac{\partial f_1^2}{\partial x^h}=0.
\]
Also, contracting (\ref{mean8}) with $v^\alpha v^\beta$ implies that
\[
\frac{\partial g^{\gamma\nu}}{\partial v^\gamma}\frac{\partial f_2^2}{\partial u^\nu}=0.
\]
Thus we have (\ref{equ}).

Conversely, let $(M_1, F_1)$ and $(M_2, F_2)$ are weakly Berwald manifolds and the relation (\ref{equ}) be hold. Then we have $E_{ij}=E_{\alpha\beta}=0$. Equation (\ref{equ}) gives us
\begin{eqnarray*}
\frac{\partial^2g^{kh}}{\partial y^j\partial
y^k}\frac{\partial f_1^2}{\partial x^h}\!\!\!\!&=&\!\!\!\!\ \frac{\partial^3g^{kh}}{\partial y^i\partial y^j\partial
y^k}\frac{\partial f_1^2}{\partial x^h}\\
\!\!\!\!&=&\!\!\!\!\ \frac{\partial^2g^{\nu\gamma}}{\partial v^\beta\partial
v^\gamma}\frac{\partial f_2^2}{\partial u^\nu}
\\
\!\!\!\!&=&\!\!\!\!\ \frac{\partial^3g^{\nu\gamma}}{\partial v^\alpha\partial v^\beta\partial
v^\gamma}\frac{\partial f_2^2}{\partial u^\nu}=0.
\end{eqnarray*}
By setting $E_{ij}=E_{\alpha\beta}=0$ and the above equation in (\ref{mean1}), (\ref{mean2}) and (\ref{mean3}),  we  can obtain $\textbf{E}_{\alpha\beta}=\textbf{E}_{ij}=\textbf{E}_{i\beta}=0$. This means that $({}_{f_2}M_1\times_{f_1}M_2, F)$ is a weakly Berwald manifold.
\end{proof}

\smallskip

Now, if $f_2$ is a constant function  on $M_2$, then (\ref{mean4}) implies that $E_{\alpha\beta}=0$. Thus we have the following.
\begin{cor}
Let $({}_{f_2}M_1\times_{f_1}M_2, F)$ be a DWP-Finsler manifold and $f_1$ be constant on $M_1$ (resp, $f_2$ be constant on $M_2$). Then $({}_{f_2}M_1\times_{f_1}M_2, F)$ is weakly Berwald if and only if $(M_1, F_1)$ and $(M_2, F_2)$ are weakly Berwald manifolds and the following holds
\[
\frac{\partial g^{\gamma\nu}}{\partial v^\gamma}\frac{\partial f_2^2}{\partial u^\nu}=0, \ \  \ (\textrm{resp},\ \  \frac{\partial g^{kh}}{\partial y^k}\frac{\partial f_1^2}{\partial x^h}=0).
\]
\end{cor}
\begin{cor}
Let $(M_1\times_{f_1}M_2, F)$ be a WP-Finsler manifold. Then $(M_1\times_{f_1}M_2, F)$ is weakly Berwald if and only if $(M_1, F_1)$ and $(M_2, F_2)$ are weakly Berwald manifolds and the following holds
\[
\frac{\partial g^{kh}}{\partial y^k}\frac{\partial f_1^2}{\partial x^h}=0.
\]
\end{cor}
\bigskip
Now, we consider DWP-Finsler manifolds with  isotropic mean Berwald  curvature. First, as a
consequence of Lemma \ref{Lemma2}, we have the following.
\begin{lem}
Every DWP-Finsler manifold $({}_{f_2}M_1\times_{f_1}M_2, F)$ has isotropic mean Berwald curvature if and only if the following hold
\begin{eqnarray}
E_{\alpha\beta}\!\!\!\!&-&\!\!\!\!\ \frac{1}{8f_1^2}\frac{\partial^3g^{\gamma\nu}}{\partial v^\beta\partial v^\alpha\partial v^\gamma}\frac{\partial f_2^2}{\partial u^\nu}F_1^2-\frac{1}{4f_2^2}g_{\alpha\beta}\frac{\partial g^{kh}}{\partial y^k}\frac{\partial f_1^2}{\partial x^h}\nonumber\\
\!\!\!\!&-&\!\!\!\!\ \frac{n+1}{2}cf_1^2F^{-1}(g_{\alpha\beta}-\frac{f_1^2}{F^2}v_\alpha v_\beta)=0,\label{wmean1}\\
E_{ij}\!\!\!\!&-&\!\!\!\!\ \frac{1}{8f_2^2}\frac{\partial^3g^{kh}}{\partial
y^i\partial y^j\partial y^k}\frac{\partial f_1^2}{\partial
x^h}F_2^2-\frac{1}{4f_1^2}g_{ij}\frac{\partial g^{\alpha\gamma}}{\partial v^\gamma}\frac{\partial f_2^2}{\partial u^\alpha}\nonumber\\
\!\!\!\!&-&\!\!\!\!\ \frac{n+1}{2}cf_2^2F^{-1}(g_{ij}-\frac{f_2^2}{F^2}y_iy_j)=0,\label{wmean2}
\end{eqnarray}
\begin{equation}
(n+1)c\frac{f_1^2f_2^2}{F^3}y_iv_\beta-\frac{1}{2f_2^2}\frac{\partial^2g^{kh}}{\partial y^k\partial
y^i}\frac{\partial f_1^2}{\partial x^h}v_\beta-\frac{1}{2f_1^2}\frac{\partial^2g^{\alpha\gamma}}{\partial v^\beta\partial
v^\gamma}\frac{\partial f_2^2}{\partial u^\alpha}y_i=0,\label{wmean3}
\end{equation}
where $c=c(\textbf{x})$ is a scalar function on $M$.
\end{lem}
\begin{thm}
Every DWP-Finsler manifold $({}_{f_2}M_1\times_{f_1}M_2, F)$  with isotropic mean Berwald curvature is a weakly Berwald manifold, provided that the following hold
\[
\frac{\partial g^{kh}}{\partial y^k}\frac{\partial f_1}{\partial x^h}=0 \ \  \  or \ \ \ \frac{\partial g^{\gamma\nu}}{\partial v^\gamma}\frac{\partial f_2}{\partial u^\nu}=0.
\]
\end{thm}
\begin{proof}
Suppose that $\frac{\partial g^{kh}}{\partial y^k}\frac{\partial f_1}{\partial x^h}=0$ and $F$ is isotropic mean Berwald DWP-Finsler metric. Then by using (\ref{wmean3}), we obtain
\[
(n+1)c\frac{f_1^2f_2^2}{F^3}v_\beta=\frac{1}{2f_1^2}\frac{\partial^2g^{\alpha\gamma}}{\partial v^\beta\partial
v^\gamma}\frac{\partial f_2^2}{\partial u^\alpha}.
\]
Differentiating the above equation with respect $y^j$ gives us
\[
\frac{(n+1)}{F^5}cf_1^2f_2^4v_\beta y_j=0.
\]
Thus, we conclude that $c=0$. This implies that $F$ is weakly Berwald metric.
\end{proof}

\section{Douglas DWP-Finsler Manifolds}
In this section, we study DWP-Finsler manifolds with vanishing Douglas curvature. We prove that every Douglas proper   DWP-Finsler manifold  is Riemannian. To prove this, we need the following.
\begin{lem}\label{Lemma4}
Let $({}_{f_2}M_1\times_{f_1}M_2, F)$ be a DWP-Finsler manifold. Then the Douglas curvature of $F$ are given as follows
\begin{eqnarray}
\textbf{D}^k_{\ ijl}\!\!\!\!&=&\!\!\!\!B^k_{\ ijl}-\frac{1}{4f_2^2}\frac{\partial^3g^{kh}}{\partial y^i\partial y^j\partial y^l}\frac{\partial f_1^2}{\partial x^h}F_2^2-\frac{2}{n+1}\Big\{E_{ij}\delta^k_l-\frac{1}{8f_2^2}\frac{\partial^3g^{sh}}{\partial
y^i\partial y^j\partial y^s}\frac{\partial f_1^2}{\partial
x^h}F_2^2\delta^k_l\nonumber\\
\!\!\!\!&&\!\!\!\!-\frac{1}{4f_1^2}g_{ij}\frac{\partial g^{\alpha\gamma}}{\partial v^\gamma}\frac{\partial f_2^2}{\partial u^\alpha}\delta^k_l+E_{il}\delta^k_j-\frac{1}{8f_2^2}\frac{\partial^3g^{sh}}{\partial
y^i\partial y^l\partial y^s}\frac{\partial f_1^2}{\partial
x^h}F_2^2\delta^k_j
\nonumber\\
\!\!\!\!&&\!\!\!\!\ -\frac{1}{4f_1^2}g_{il}\frac{\partial g^{\alpha\gamma}}{\partial v^\gamma}\frac{\partial f_2^2}{\partial u^\alpha}\delta^k_j+E_{jl}\delta^k_i-\frac{1}{8f_2^2}\frac{\partial^3g^{sh}}{\partial
y^j\partial y^l\partial y^s}\frac{\partial f_1^2}{\partial
x^h}F_2^2\delta^k_i\nonumber\\
\!\!\!\!&&\!\!\!\!\
-\frac{1}{4f_1^2}g_{jl}\frac{\partial g^{\alpha\gamma}}{\partial v^\gamma}\frac{\partial f_2^2}{\partial u^\alpha}\delta^k_i-\frac{1}{4f_1^2}\frac{\partial g_{ij}}{\partial y^l}\frac{\partial g^{\alpha\gamma}}{\partial v^\gamma}\frac{\partial f_2^2}{\partial u^\alpha}y^k \nonumber\\
\!\!\!\!&&\!\!\!\!-\frac{1}{8f_2^2}\frac{\partial^4g^{sh}}{\partial y^l\partial
y^i\partial y^j\partial y^s}\frac{\partial f_1^2}{\partial
x^h}F_2^2y^k  +\frac{\partial E_{ij}}{\partial y^l}y^k    \Big\},
\end{eqnarray}
\begin{eqnarray}
\textbf{D}^k_{\ i\beta l}\!\!\!\!&=&\!\!\!\!-\frac{1}{4f_2^2}\frac{\partial^2g^{kh}}{\partial y^l\partial y^i}\frac{\partial f_1^2}{\partial x^h}\frac{\partial F_2^2}{\partial v^\beta}+\frac{2}{n+1}\Big\{\frac{1}{4f_2^2}\frac{\partial^2g^{sh}}{\partial y^s\partial
y^i}\frac{\partial f_1^2}{\partial x^h}\delta^k_lv_\beta\nonumber\\
\!\!\!\!&&\!\!\!\!+\frac{1}{4f_1^2}\frac{\partial^2g^{\alpha\gamma}}{\partial v^\beta\partial
v^\gamma}\frac{\partial f_2^2}{\partial u^\alpha}\delta^k_ly_i+\frac{1}{4f_2^2}\frac{\partial^2g^{sh}}{\partial y^s\partial
y^l}\frac{\partial f_1^2}{\partial x^h}\delta^k_iv_\beta\nonumber\\
\!\!\!\!&&\!\!\!\!+\frac{1}{4f_1^2}\frac{\partial^2g^{\alpha\gamma}}{\partial v^\beta\partial
v^\gamma}\frac{\partial f_2^2}{\partial u^\alpha}\delta^k_iy_l
+\frac{1}{4f_2^2}\frac{\partial^3g^{sh}}{\partial y^l\partial y^s\partial
y^i}\frac{\partial f_1^2}{\partial x^h}y^kv_\beta\nonumber\\
\!\!\!\!&&\!\!\!\!+\frac{1}{4f_1^2}\frac{\partial^2g^{\alpha\gamma}}{\partial v^\beta\partial
v^\gamma}\frac{\partial f_2^2}{\partial u^\alpha}y^kg_{il}\Big\},
\end{eqnarray}
\begin{eqnarray}
\textbf{D}^k_{\ \alpha\beta l}\!\!\!\!&=&\!\!\!\!-\frac{1}{2f_2^2}g_{\alpha\beta}\frac{\partial g^{kh}}{\partial y^l}\frac{\partial f_1^2}{\partial x^h}-\frac{2}{n+1}\Big\{E_{\alpha\beta}\delta^k_l-\frac{1}{4f_2^2}g_{\alpha\beta}\frac{\partial g^{sh}}{\partial y^s}\frac{\partial f_1^2}{\partial x^h}\delta^k_l\nonumber\\
\!\!\!\!&&\!\!\!\!-\frac{1}{8f_1^2}\frac{\partial^3g^{\gamma\nu}}{\partial v^\beta\partial v^\alpha\partial v^\gamma}\frac{\partial f_2^2}{\partial u^\nu}F_1^2\delta^k_l-\frac{1}{4f_1^2}\frac{\partial^3g^{\gamma\nu}}{\partial v^\beta\partial v^\alpha\partial v^\gamma}\frac{\partial f_2^2}{\partial u^\nu}y_ly^k\nonumber\\
\!\!\!\!&&\!\!\!\!-\frac{1}{4f_2^2}g_{\alpha\beta}\frac{\partial^2 g^{sh}}{\partial y^l\partial y^s}\frac{\partial f_1^2}{\partial x^h}y^k
\Big\},
\end{eqnarray}
\begin{eqnarray}
\textbf{D}^k_{\ \alpha\beta\lambda}=-\frac{1}{f_2^2}C_{\alpha\beta\lambda}g^{kh}\frac{\partial f_1^2}{\partial x^h} \!\!\!\!&-&\!\!\!\! \frac{2}{n+1}\Big\{\frac{\partial E_{\alpha\beta}}{\partial v^\lambda}y^k-\frac{1}{8f_1^2}\frac{\partial^4g^{\gamma\nu}}{\partial v^\lambda\partial v^\beta\partial v^\alpha\partial v^\gamma}\frac{\partial f_2^2}{\partial u^\nu}F_1^2y^k\nonumber\\
\!\!\!\!&-&\!\!\!\!\frac{1}{4f_2^2}\frac{\partial g_{\alpha\beta}}{\partial v^\lambda}\frac{\partial g^{sh}}{\partial y^s}\frac{\partial f_1^2}{\partial x^h}y^k
\Big\},\label{doug4}
\end{eqnarray}
\begin{eqnarray}
\textbf{D}^\gamma_{\ \alpha\beta\lambda}= B^\gamma_{\alpha\beta\lambda}\!\!\!\!&-&\!\!\!\!\  \frac{1}{4f_1^2}
\frac{\partial^3g^{\gamma\nu}}{\partial v^\beta\partial v^\alpha\partial v^\lambda}\frac{\partial f_2^2}{\partial u^\nu}F_1^2+\frac{2}{n+1}\Big\{\frac{1}{8f_1^2}\frac{\partial^3g^{\mu\nu}}{\partial v^\beta\partial v^\alpha\partial v^\mu}\frac{\partial f_2^2}{\partial u^\nu}F_1^2\delta^\gamma_\lambda\nonumber\\
\!\!\!\!&-&\!\!\!\!\ E_{\alpha\beta}\delta^\gamma_\lambda+\frac{1}{4f_2^2}g_{\alpha\beta}\frac{\partial g^{kh}}{\partial y^k}\frac{\partial f_1^2}{\partial x^h}\delta^\gamma_\lambda +\frac{1}{8f_1^2}\frac{\partial^3g^{\mu\nu}}{\partial v^\lambda\partial v^\alpha\partial v^\mu}\frac{\partial f_2^2}{\partial u^\nu}F_1^2\delta^\gamma_\beta
\nonumber\\
\!\!\!\!&-&\!\!\!\!\  E_{\alpha\lambda}\delta^\gamma_\beta+\frac{1}{4f_2^2}g_{\alpha\lambda}\frac{\partial g^{kh}}{\partial y^k}\frac{\partial f_1^2}{\partial x^h}\delta^\gamma_\beta+\frac{1}{8f_1^2}\frac{\partial^3g^{\mu\nu}}{\partial v^\lambda\partial v^\beta\partial v^\mu}\frac{\partial f_2^2}{\partial u^\nu}F_1^2\delta^\gamma_\alpha
\nonumber\\
\!\!\!\!&-&\!\!\!\!\ E_{\beta\lambda}\delta^\gamma_\alpha+\frac{1}{4f_2^2}g_{\beta\lambda}\frac{\partial g^{kh}}{\partial y^k}\frac{\partial f_1^2}{\partial x^h}\delta^\gamma_\alpha+\frac{1}{8f_1^2}\frac{\partial^4g^{\mu\nu}}{\partial v^\lambda\partial v^\beta\partial v^\alpha\partial v^\mu}\frac{\partial f_2^2}{\partial u^\nu}F_1^2v^\gamma\nonumber\\
\!\!\!\!&+&\!\!\!\!\ \frac{1}{4f_2^2}\frac{g_{\alpha\beta}}{\partial v^\gamma}\frac{\partial g^{kh}}{\partial y^k}\frac{\partial f_1^2}{\partial x^h}v^\gamma-\frac{\partial E_{\alpha\beta}}{\partial v^\lambda}v^\gamma\Big\},
\end{eqnarray}
\begin{eqnarray}
\textbf{D}^\gamma_{\ i\beta\lambda}\!\!\!\!&=&\!\!\!\!-\frac{1}{4f_1^2}\frac{\partial^2g^{\alpha\gamma}}{\partial v^\beta\partial v^\lambda}\frac{\partial f_2^2}{\partial u^\alpha}\frac{\partial F_1^2}{\partial y^i}+\frac{2}{n+1}\Big\{\frac{1}{4f_2^2}\frac{\partial^2g^{kh}}{\partial y^k\partial
y^i}\frac{\partial f_1^2}{\partial x^h}\delta^\gamma_\lambda v_\beta\nonumber\\
\!\!\!\!&&\!\!\!\!
+\frac{1}{4f_1^2}\frac{\partial^2g^{\alpha\mu}}{\partial v^\beta\partial
v^\mu}\frac{\partial f_2^2}{\partial u^\alpha}\delta^\gamma_\lambda y_i+\frac{1}{4f_2^2}\frac{\partial^2g^{kh}}{\partial y^k\partial
y^i}\frac{\partial f_1^2}{\partial x^h}\delta^\gamma_\beta v_\lambda
+\frac{1}{4f_1^2}\frac{\partial^2g^{\alpha\mu}}{\partial v^\lambda\partial
v^\mu}\frac{\partial f_2^2}{\partial u^\alpha}\delta^\gamma_\beta y_i\nonumber\\
\!\!\!\!&&\!\!\!\!+\frac{1}{4f_2^2}\frac{\partial^2g^{kh}}{\partial y^k\partial
y^i}\frac{\partial f_1^2}{\partial x^h}g_{\beta\lambda}v^\gamma
+\frac{1}{4f_1^2}\frac{\partial^3g^{\alpha\mu}}{\partial v^\lambda\partial v^\beta\partial
v^\mu}\frac{\partial f_2^2}{\partial u^\alpha} v^\gamma y_i
\Big\},
\end{eqnarray}
\begin{eqnarray}
\textbf{D}^\gamma_{\ ij\lambda}\!\!\!\!&=&\!\!\!\!-\frac{1}{2f_1^2}g_{ij}\frac{\partial g^{\alpha\gamma}}{\partial v^\lambda}\frac{\partial f_2^2}{\partial u^\alpha}-\frac{2}{n+1}\Big\{E_{ij}\delta^\gamma_\lambda-\frac{1}{8f_2^2}\frac{\partial^3g^{kh}}{\partial
y^i\partial y^j\partial y^k}\frac{\partial f_1^2}{\partial
x^h}F_2^2\delta^\gamma_\lambda\nonumber\\
\!\!\!\!&&\!\!\!\!-\frac{1}{4f_1^2}g_{ij}\frac{\partial g^{\alpha\mu}}{\partial v^\mu}\frac{\partial f_2^2}{\partial u^\alpha}\delta^\gamma_\lambda-\frac{1}{8f_2^2}\frac{\partial^3g^{kh}}{\partial
y^i\partial y^j\partial y^k}\frac{\partial f_1^2}{\partial
x^h}\frac{\partial F_2^2}{\partial v^\lambda}v^\gamma\nonumber\\
\!\!\!\!&&\!\!\!\!-\frac{1}{4f_1^2}g_{ij}\frac{\partial^2 g^{\alpha\mu}}{\partial v^\lambda\partial v^\mu}\frac{\partial f_2^2}{\partial u^\alpha}v^\gamma
\Big\},
\end{eqnarray}
\begin{eqnarray}
\textbf{D}^\gamma_{\ ijk}=-\frac{1}{f_1^2}C_{ijk}g^{\alpha\gamma}\frac{\partial f_2^2}{\partial u^\alpha}\!\!\!\!&-&\!\!\!\! \frac{2}{n+1}\Big\{\frac{\partial E_{ij}}{\partial y^k}v^\gamma-\frac{1}{8f_2^2}\frac{\partial^4g^{sh}}{\partial y^k\partial
y^i\partial y^j\partial y^s}\frac{\partial f_1^2}{\partial
x^h}F_2^2v^\gamma\nonumber\\
\!\!\!\!&-&\!\!\!\! \frac{1}{4f_1^2}\frac{\partial g_{ij}}{\partial y^k}\frac{\partial g^{\alpha\mu}}{\partial v^\mu}\frac{\partial f_2^2}{\partial u^\alpha}v^\gamma\Big\}.\label{doug8}
\end{eqnarray}
\end{lem}
\begin{proof}
By some lengthy calculations and using  Lemmas \ref{Lemma1} and \ref{Lemma2},  we can  get the proof.
\end{proof}

\bigskip

\begin{thm}\label{thm4}
Every  proper DWP-Finsler manifold $({}_{f_2}M_1\times{}_{f_1}M_2,F)$ with vanishing Douglas  curvature  is a Riemannian manifold.
\end{thm}
\begin{proof}
Suppose that the Douglas curvature of $({}_{f_2}M_1\times{}_{f_1}M_2,F)$  vanishes, i.e., $\textbf{D}^d_{abc}=0$. Then by contracting (\ref{doug8}) with $y^k$ we obtain
\begin{equation}
E_{ij}=\frac{3}{8f_2^2}\frac{\partial^3g^{sh}}{\partial
y^i\partial y^j\partial y^s}\frac{\partial f_1^2}{\partial
x^h}F_2^2.\label{doug9}
\end{equation}
Since $E_{ij}$ is a function of $(x,y)$, then by differentiating the above equation with respect to $v^{\alpha}$, we get
\begin{equation}
\frac{\partial^3g^{sh}}{\partial
y^i\partial y^j\partial y^s}\frac{\partial f_1^2}{\partial
x^h}=0.\label{doug10}
\end{equation}
Putting the above equation into (\ref{doug9}) gives us $E_{ij}=0$. Further, (\ref{doug10}) implies that
 \begin{equation}
\frac{\partial^4g^{sh}}{\partial y^k\partial
y^i\partial y^j\partial y^s}\frac{\partial f_1^2}{\partial
x^h}=\frac{\partial^2g^{sh}}{\partial y^j\partial y^s}\frac{\partial f_1^2}{\partial
x^h}=\frac{\partial g^{sh}}{\partial y^s}\frac{\partial f_1^2}{\partial
x^h}=0.\label{doug11}
\end{equation}
In the similar way, we can conclude that $E_{\alpha\beta}=0$ and
\begin{eqnarray}
\nonumber \frac{\partial^4g^{\gamma\nu}}{\partial v^\lambda\partial v^\beta\partial v^\alpha\partial v^\gamma}\frac{\partial f_2^2}{\partial u^\nu}\!\!\!\!&=&\!\!\!\!\ \frac{\partial^3g^{\gamma\nu}}{\partial v^\beta\partial v^\alpha\partial v^\gamma}\frac{\partial f_2^2}{\partial u^\nu}\\
\!\!\!\!&=&\!\!\!\!\ \frac{\partial^2g^{\gamma\nu}}{\partial v^\alpha\partial v^\gamma}\frac{\partial f_2^2}{\partial u^\nu}=\frac{\partial g^{\gamma\nu}}{\partial v^\gamma}\frac{\partial f_2^2}{\partial u^\nu}=0.\label{doug12}
\end{eqnarray}
Setting (\ref{doug10})-(\ref{doug12}) and $E_{ij}=E_{\alpha\beta}=0$ into (\ref{doug4}) and (\ref{doug8}) imply that $C_{ijk}=C_{\alpha\beta\lambda}=0$. Therefore $(M_1, F_1)$ and $(M_2, F_2)$ are Riemannian and consequently, $({}_{f_2}M_1\times{}_{f_1}M_2,F)$ is Riemannian.
\end{proof}

\smallskip
Using the Theorem \ref{thm4}, we have the following.
\begin{cor}
Let $({}_{f_2}M_1\times{}_{f_1}M_2,F)$ be a DWP-Finsler manifold.\\
(i) If $f_2$ is constant on $M_2$. Then $F$ is Douglas metric if and only if $F_2$ is Riemannian metric, $F_1$ is Berwald metric and $\frac{\partial g^{sh}}{\partial y^s}\frac{\partial f_1^2}{\partial
x^h}=0$.\\
(ii) If $f_1$ is constant on $M_1$. Then $F$ is Douglas metric if and only if $F_1$ is Riemannian metric, $F_2$ is Berwald metric and $\frac{\partial g^{\gamma\lambda}}{\partial v^\gamma}\frac{\partial f_2^2}{\partial
u^\lambda}=0$.
\end{cor}

\smallskip
Finally, we consider warped product Finsler manifold with vanishing Douglas curvature and obtain the following.

\begin{cor}
The WP-Finsler manifold $(M_1\times_{f_1}M_2, F)$ is Douglas
manifold if and only if $F_2$ is Riemannian metric, $F_1$ is
Berwald metric and $\frac{\partial g^{sh}}{\partial
y^s}\frac{\partial f_1^2}{\partial x^h}=0$.
\end{cor}
\begin{proof}
By the Lemma \ref{Lemma4}, we can get the proof.
\end{proof}

\section{Relatively Isotropic Landsberg DWP-Finsler Manifolds}
In this section, we prove that on a  proper DWP-Finsler manifold  the notions being Landsberg manifold and being Berwald manifold are equivalent. Then we study DWP-Finsler metrics with relatively isotropic Landsberg curvature.
\begin{lem}\label{Lemma5}
Let $({}_{f_2}M_1\times_{f_1}M_2, F)$ be a DWP-Finsler manifold. Then the Landsberg curvature of $F$ are given as follows
\begin{eqnarray}
\textbf{L}_{ijk}\!\!\!\!&=&\!\!\!\!\  f_2^2L_{ijk}+\frac{1}{8}f_2^2y_l\frac{\partial^3g^{lh}}{\partial
y^i\partial y^j\partial y^k}\frac{\partial f_1^2}{\partial
x^h}F_2^2+\frac{1}{2}C_{ijk}v^{\alpha}\frac{\partial
f_2^2}{\partial u^\alpha},\label{lan19}
\\
\textbf{L}_{ij\lambda}\!\!\!\!&=&\!\!\!\!\ \frac{1}{4}y_l\frac{\partial^2g^{lh}}{\partial
y^i\partial y^j}\frac{\partial f_1^2}{\partial
x^h}v_{\lambda}+\frac{1}{4}g_{ij}v_{\gamma}\frac{\partial
g^{\alpha\gamma}}{\partial v^\lambda}\frac{\partial
f_2^2}{\partial u^\alpha},
\\
\textbf{L}_{i\beta\lambda}\!\!\!\!&=&\!\!\!\!\ \frac{1}{4}y_l\frac{\partial
g^{lh}}{\partial y^i}\frac{\partial f_1^2}{\partial
x^h}g_{\beta\lambda}+\frac{1}{4}v_{\gamma}\frac{\partial^2
g^{\alpha\gamma}}{\partial v^\beta\partial
v^\lambda}\frac{\partial f_2^2}{\partial u^\alpha}y_i,
\\
\textbf{L}_{\alpha\beta\lambda}\!\!\!\!&=&\!\!\!\!\ f_1^2L_{\alpha\beta\lambda}+\frac{1}{8}f_1^2v_\gamma\frac{\partial^3g^{\gamma\nu}}{\partial
v^\alpha\partial v^\beta\partial v^\lambda}\frac{\partial
f_2^2}{\partial
u^\nu}F_1^2+\frac{1}{2}C_{\alpha\beta\lambda}y^h\frac{\partial
f_1^2}{\partial x^h}.
\end{eqnarray}
\end{lem}
\begin{proof}
By the definition of Landsberg curvature and (\ref{feri})-(\ref{2a}), we can get the proof.
\end{proof}

\smallskip

\begin{prop}\label{prop6}
Every  proper DWP-Finsler manifold
$({}_{f_2}M_1\times{}_{f_1}M_2,F)$ with vanishing Landsberg
curvature is  Riemannian.
\end{prop}
\begin{proof}
Let the Landsberg curvature tensor of
$({}_{f_2}M_1\times{}_{f_1}M_2,F)$ be zero. Then by using
(\ref{lan19}) we obtain
\begin{equation}
f_2^2L_{ijk}+\frac{1}{8}f_2^2y_l\frac{\partial^3g^{lh}}{\partial
y^i\partial y^j\partial y^k}\frac{\partial f_1^2}{\partial
x^h}F_2^2+\frac{1}{2}C_{ijk}v^{\alpha}\frac{\partial
f_2^2}{\partial u^\alpha}=0.\label{lan5}
\end{equation}
Differentiating (\ref{lan5}) with respect to $v^\gamma$ implies that
\begin{equation}
\frac{1}{4}f_2^2y_l\frac{\partial^3g^{lh}}{\partial y^i\partial
y^j\partial y^k}\frac{\partial f_1^2}{\partial
x^h}v_\gamma+\frac{1}{2}C_{ijk}\frac{\partial f_2^2}{\partial
u^\gamma}=0.\label{lan5.1}
\end{equation}
By differentiating (\ref{lan5.1}) with respect $v^\lambda$, we have
\[
y_l\frac{\partial^3g^{lh}}{\partial y^i\partial y^j\partial
y^k}\frac{\partial f_1^2}{\partial x^h}g_{\gamma\lambda}=0,
\]
and consequently
\[
y_l\frac{\partial^3g^{lh}}{\partial y^i\partial
y^j\partial y^k}\frac{\partial f_1^2}{\partial x^h}=0.
\]
Thus (\ref{lan5}) reduces to following
\begin{equation}
f_2^2L_{ijk}+\frac{1}{2}C_{ijk}v^{\alpha}\frac{\partial
f_2^2}{\partial u^\alpha}=0.\label{lan6}
\end{equation}
Differentiating (\ref{lan6}) with respect to $v^\beta$ gives us
\[
C_{ijk}\frac{\partial f_2^2}{\partial u^\beta}=0,
\]
and consequently $C_{ijk}=0$. Thus $(M_1, F_1)$ is a Riemannian manifold. By the similar way, we can conclude that $(M_2, F_2)$ is
Riemannian.
\end{proof}
Using Proposition \ref{prop6} and Lemma \ref{Lemma1}, we get the following.
\begin{thm} \label{Lan-Ber}
A proper DWP-Finsler manifold is Landsbergin if and only if it is Berwaldian.
\end{thm}
Now, let $f_1$ be non-constant on $M_1$ and $f_2$ be constant on
$M_2$. Then by the  similar way used in to the proof of Proposition \ref{prop6}, we conclude that
$(M_2, F_2)$ is a Riemannain manifold. Also, from (\ref{lan6}) we conclude $L_{ijk}=0$, because $f_2$ is constant. Thus we have the following.
\begin{thm}\label{thm7}
Let $({}_{f_2}M_1\times{}_{f_1}M_2,F)$ be a proper DWP-Finsler manifold. Then\\
(i)\ \ If $f_2$ is constant and $f_1$ is not constant, then
$({}_{f_2}M_1\times{}_{f_1}M_2,F)$ is a Landsberg manifold if and
only if $(M_1, F_1)$ is  a Landsberg manifold, $(M_2, F_2)$ is Riemannian and the following  holds
\be
y_l\frac{\partial^3g^{lh}}{\partial y^i\partial y^j\partial
y^k}\frac{\partial f_1}{\partial x^h}=0.\label{Lan}
\ee
(ii)\ \ If $f_1$ is constant and $f_2$ is not constant, then
$({}_{f_2}M_1\times{}_{f_1}M_2,F)$ is a Landsberg manifold if and
only if $(M_2, F_2)$ is a Landsberg manifold, $(M_1, F_1)$ is Riemannian and the following is hold
\be
v_\gamma\frac{\partial^3g^{\gamma\nu}}{\partial v^\alpha\partial
v^\beta\partial v^\lambda}\frac{\partial f_2}{\partial
u^\nu}=0.
\ee
\end{thm}

\bigskip
By Theorem \ref{thm7}, we have the following.
\begin{cor}
The WP-Finsler manifold $(M_1\times{}_{f_1}M_2,F)$ is
Landsberg manifold if and only if $(M_1, F_1)$ is Landsberg,
$(M_2, F_2)$ is Riemannian and
\be
C^h_{\,\,kj}\frac{\partial f_1}{\partial x^h}=0,\label{Lan33}
\ee
\end{cor}
\begin{proof}
It suffices to show that (\ref{Lan}) implies (\ref{Lan33}).  Multiplying  (\ref{Lan}) with $y^i$ implies that
\be
y_l\frac{\partial C^{lh}_{\,\,\,\,\,j}}{\partial
y^k}\frac{\partial f_1}{\partial x^h}=0,\label{Lan22}
\ee
Using $y_lC^{lh}_{\,\,\,\,j}=0$ and $\frac{\partial y_l}{\partial y^k}=g_{lk}$, one can obtain (\ref{Lan33}).
\end{proof}

\bigskip

Now, we are going to consider DWP-Finsler manifold with relatively isotropic Landsberg metric.

\begin{thm}\label{thm9}
Let $({}_{f_2}M_1\times{}_{f_1}M_2,F)$ be a DWP-Finsler manifold. Suppose that $F$ is a relatively isotropic Landsberg metric. Then $F$ is a  Landsberg metric.
\end{thm}
\begin{proof}
Let $({}_{f_2}M_1\times{}_{f_1}M_2,F)$ be a relatively isotropic
Landsberg manifold. Then by (\ref{lan19}), we have
\be
f_2^2L_{ijk}+\frac{1}{8}f_2^2y_l\frac{\partial^3g^{lh}}{\partial
y^i\partial y^j\partial y^k}\frac{\partial f_1^2}{\partial
x^h}F_2^2+\frac{1}{2}C_{ijk}v^{\alpha}\frac{\partial
f_2^2}{\partial u^\alpha}=cFf_2^2C_{ijk}.\label{ril1}
\ee
By differentiating (\ref{ril1}) with respect $v^\gamma$ and $v^\lambda$, one can obtain the following
\be
\frac{1}{4}y_l\frac{\partial^3g^{lh}}{\partial y^i\partial
y^j\partial y^k}\frac{\partial f_1^2}{\partial
x^h}g_{\gamma\lambda}=cF^{-1}{\bf{h}}_{\gamma\lambda}C_{ijk}.\label{ril2}
\ee
Contracting (\ref{ril2}) with $g^{\gamma\lambda}$ implies that
\be
\frac{n_2}{4}y_l\frac{\partial^3g^{lh}}{\partial y^i\partial
y^j\partial y^k}\frac{\partial f_1^2}{\partial
x^h}=(n-1)cF^{-1}C_{ijk}.\label{ril3}
\ee
Differentiating (\ref{ril3}) with respect to $v^\beta$ gives us
\[
(n-1)c(F^{-1})_{v^{\beta}}C_{ijk}=0,
\]
and then  $c=0$. Thus $F$ reduces to a  Landsberg metric.
\end{proof}

By Proposition \ref{prop6} and Theorem \ref{thm9}, we conclude the following
\begin{cor}\label{cor5}
Every proper DWP-Finsler manifold with relatively isotropic Landsberg curvature is Riemannian.
\end{cor}

\section{Relatively Isotropic Mean Landsberg DWP-Finsler Manifolds}
In this section, we consider DWP-Finsler metrics with relatively isotropic mean Landsberg curvature. First, by the definition of mean  Landsberg curvature and  Lemma \ref{Lemma5}, we can get the following.
\begin{lem}
Let $({}_{f_2}M_1\times_{f_1}M_2, F)$ be a DWP-Finsler manifold. Then the mean Landsberg curvature of $F$ are given as follows
\begin{eqnarray}
\nonumber \textbf{J}_i\!\!\!\!&=&\!\!\!\!\!\!\!\ \frac{1}{f_2^2}g^{jk}\textbf{L}_{ijk}+\frac{1}{f_1^2}g^{\beta\lambda}\textbf{L}_{i\beta\lambda}\\
\!\!\!\!&=&\!\!\!\!\!\!\!\!\ J_i+\frac{y_lg^{jk}}{8}\frac{\partial^3g^{lh}}{\partial
y^i\partial y^j\partial y^k}\frac{\partial f_1^2}{\partial
x^h}F_2^2+\frac{I_iv^{\nu}}{2f_2^2}\frac{\partial
f_2^2}{\partial u^\nu}+\frac{v_{\gamma}g^{\beta\lambda}}{4}\frac{\partial^2
g^{\nu\gamma}}{\partial v^\beta\partial
v^\lambda}\frac{\partial f_2^2}{\partial
u^\nu}y_i,\label{weakl1}
\\
\nonumber \textbf{J}_\alpha \!\!\!\!&=&\!\!\!\!\!\!\!\!\!\ \frac{1}{f_2^2}g^{jk}\textbf{L}_{\alpha
jk}+\frac{1}{f_1^2}g^{\beta\lambda}\textbf{L}_{\alpha\beta\lambda}\\
\!\!\!\!&=&\!\!\!\!\!\!\!\!\!\ J_\alpha+\frac{v_\gamma g^{\beta\lambda}}{8}\frac{\partial^3g^{\gamma
\nu}}{\partial v^\alpha\partial v^\beta\partial
v^\lambda}\frac{\partial f_2^2}{\partial
u^\nu}F_1^2+\frac{I_\alpha y^h}{2f_1^2}\frac{\partial
f_1^2}{\partial x^h}+\frac{y_lg^{jk}}{4}\frac{\partial^2 g^{lh}}{\partial
y^j\partial y^k}\frac{\partial f_1^2}{\partial
x^h}v_\alpha.\label{weakl2}
\end{eqnarray}
\end{lem}
\bigskip

\begin{thm}\label{THM11}
Let $({}_{f_2}M_1\times{}_{f_1}M_2,F)$ be a DWP-Finsler manifold. Then\\
(i) If  $f_1$ is constant and $f_2$ is not constant, then $({}_{f_2}M_1\times{}_{f_1}M_2,F)$ is a
weakly Landsberg manifold if and only if $(M_1, F_1)$ is
Riemannian, $(M_2, F_2)$ is weakly Landsbergian and  the following  holds
\[
v_\gamma\frac{\partial^3g^{\gamma \nu}}{\partial
v^\alpha\partial v^\beta\partial
v^\lambda}g^{\beta\lambda}\frac{\partial f_2^2}{\partial u^\nu}=0.
\]
(ii) If  $f_2$ is constant and $f_1$ is not constant,  then $({}_{f_2}M_1\times{}_{f_1}M_2,F)$ is a
weakly Landsberg manifold if and only if $(M_1, F_1)$ is weakly Landsbergian, $(M_2, F_2)$ is Riemannian and  the following  holds
\[
y_l\frac{\partial^3g^{lh}}{\partial y^i\partial y^j\partial y^k}g^{jk}\frac{\partial f_1^2}{\partial x^h}=0.
\]
\end{thm}
\begin{proof}
Let $({}_{f_2}M_1\times{}_{f_1}M_2,F)$ be a weakly Landsberg manifold and $f_1$ be  constant on $M_1$. Then by (\ref{weakl1}) and (\ref{weakl2}),  we have
\begin{eqnarray}
&&J_i+\frac{1}{2f_2^2}I_iv^{\nu}\frac{\partial f_2^2}{\partial
u^\nu}+\frac{1}{4}v_{\gamma}\frac{\partial^2
g^{\nu\gamma}}{\partial v^\beta\partial
v^\lambda}g^{\beta\lambda}\frac{\partial f_2^2}{\partial
u^\nu}y_i=0,\label{weakl3}
\\
&&J_\alpha+\frac{1}{8}v_\gamma\frac{\partial^3g^{\gamma
\nu}}{\partial v^\alpha\partial v^\beta\partial
v^\lambda}g^{\beta\lambda}\frac{\partial f_2^2}{\partial
u^\nu}F_1^2=0.\label{weakl4}
\end{eqnarray}
By differentiating (\ref{weakl4}) with respect to $y^i$, we get
\begin{equation}
v_\gamma\frac{\partial^3g^{\gamma \nu}}{\partial v^\alpha\partial
v^\beta\partial v^\lambda}g^{\beta\lambda}\frac{\partial
f_2^2}{\partial u^\nu}=0.\label{weakl5}
\end{equation}
Contracting (\ref{weakl5}) with $v^\alpha$ gives us
\be
v_\gamma\frac{\partial^2g^{\gamma \nu}}{\partial v^\beta\partial
v^\lambda}g^{\beta\lambda}\frac{\partial f_2^2}{\partial u^\nu}=0.\label{weakl6}
\ee
Setting  (\ref{weakl6}) into (\ref{weakl3}) implies that
\be
J_i+\frac{1}{2f_2^2}I_iv^{\nu}\frac{\partial f_2^2}{\partial u^\nu}=0.\label{weakl7}
\ee
By differentiating   (\ref{weakl7}) with respect $v^\beta$,  we conclude that $I_i=0$, i.e., $(M_1, F_1)$ is a Riemannian manifold. By setting
(\ref{weakl5}) into (\ref{weakl2}),  we get  $J_\alpha=0$, i.e.,
$(M_2, F_2)$ is weakly Landsberg manifold. The converse of the
assertion is obvious.
\end{proof}

\smallskip

By Theorem \ref{THM11}, we can conclude the following.
\begin{cor}
The proper WP-Finsler manifold
$(M_1\times{}_{f_1}M_2,F)$ is weakly Landsberg manifold if and
only if $(M_1, F_1)$ is weakly Landsberg, $(M_2, F_2)$ is
Riemannian and the following  holds
\[
y_l\frac{\partial^3g^{lh}}{\partial y^i\partial
y^j\partial y^k}g^{jk}\frac{\partial f_1^2}{\partial x^h}=0.
\]
\end{cor}

\bigskip
Now, we are going to consider DWP-Finsler manifold with relatively isotropic mean Landsberg curvature.

\begin{thm}
Let $({}_{f_2}M_1\times{}_{f_1}M_2,F)$ be a DWP-Finsler manifold. If $f_1$ is constant on $M_1$ (resp, $f_2$ is constant on $M_2$), then the  DWP-Finsler manifold with relatively isotropic mean Landsberg curvature  is a weakly Landsberg manifold.
\end{thm}
\begin{proof}
Let $({}_{f_2}M_1\times{}_{f_1}M_2,F)$ be a relatively isotropic
mean Landsberg manifold  and $f_1$ be a constant on $M_1$. Then by (\ref{weakl4}) we have
\begin{equation}
J_\alpha+\frac{1}{8}v_\gamma\frac{\partial^3g^{\gamma
\nu}}{\partial v^\alpha\partial v^\beta\partial
v^\lambda}g^{\beta\lambda}\frac{\partial f_2^2}{\partial
u^\nu}F_1^2=cFI_\alpha.\label{ref}
\end{equation}
Differentiating (\ref{ref}) with respect to $y^k$  implies that
\be
\frac{1}{4}v_\gamma\frac{\partial^3g^{\gamma \nu}}{\partial
v^\alpha\partial v^\beta\partial
v^\lambda}g^{\beta\lambda}\frac{\partial f_2^2}{\partial
u^\nu}y_k=c\frac{f_2^2y_k}{F}I_\alpha.\label{ref2}
\ee
Contracting (\ref{ref2}) with $y^k$ gives us
\be
\frac{1}{4}v_\gamma\frac{\partial^3g^{\gamma \nu}}{\partial
v^\alpha\partial v^\beta\partial
v^\lambda}g^{\beta\lambda}\frac{\partial f_2^2}{\partial
u^\nu}F_1^2=c\frac{f_2^2F_1^2}{F}I_\alpha.\label{ref3}
\ee
By setting  (\ref{ref3}) in (\ref{ref}), it follows that
\be
2J_\alpha+c(\frac{f_2^2F_1^2}{F}-2F)I_\alpha=0.\label{ref4}
\ee
By differentiating  (\ref{ref4})  with respect to $y^k$, we obtain $c\frac{f_2^4F_1^2}{F^3}y_kI_\alpha=0$. Therefore, $c=0$ and $F$ reduces to a weakly Landsberg metric.
\end{proof}
\begin{cor}
Every WP-manifold $(M_1\times{}_{f_1}M_2,F)$ with relatively isotropic mean Landsberg curvature is a weakly Landsberg manifold.
\end{cor}

\section{Locally Dually Flat DWP-Finsler Manifolds}
Dually flat Finsler metrics form a special and valuable class of Finsler metrics in Finsler information geometry, which play a very   important role in studying flat Finsler information structure. In this section, we study locally dually flat  DWP-Finsler metrics. It is remarkable that, a Finsler metric $F=F(x,y)$ on a manifold $M$ is said to be locally dually flat if at any point there is a standard coordinate system
 $(x^i,y^i)$ in $TM$ such that it  satisfies
\be
\frac{\partial^2 F^2}{\partial x^k \partial y^l}\ y^k=2\ \frac{\partial F^2}{\partial x^l}.\label{1}
\ee
In this case, the coordinate $(x^i)$ is called an adapted local coordinate system.

\begin{thm}\label{THM1}
Let $({}_{f_2}M_1\times_{f_1}M_2, F)$ be a DWP-Finsler manifold. Then $F$ is locally
dually flat if and only if $F_1$ and $F_2$ are locally
dually flat and $f_1$ and $f_2$ are constant.
\end{thm}
\begin{proof}
Let $({}_{f_1}M_1\times{}_{f_2}M_2,F)$ be a locally dually flat doubly DWP-Finsler manifold. Then we have
\begin{equation}
f_2^2\frac{\partial^2F_1^2}{\partial x^k\partial
y^l}y^k+\frac{\partial f_2^2}{\partial u^\alpha}\frac{\partial
F_1^2}{\partial y^l}v^\alpha=2f_2^2\frac{\partial F_1^2}{\partial x^l}+2\frac{\partial
f_1^2}{\partial x^l}F_2^2,\label{Im1}
\end{equation}
\begin{equation}
\frac{\partial f_1^2}{\partial x^k}\frac{\partial F_2^2}{\partial
v^\beta}y^k+f_1^2\frac{\partial^2F_2^2}{\partial u^\alpha\partial
v^\beta}v^\alpha=2f_1^2\frac{\partial F_2^2}{\partial
u^\beta}+2\frac{\partial f_2^2}{\partial u^\beta}F_1^2. \label{Im2}
\end{equation}
Taking derivative with respect to $v^{\gamma}$ and then to $y^k$  from (\ref{Im1}) and using
non-singularity of $g_{ij}$ yield
\[
\frac{\partial f_2}{\partial u^{\gamma}}=0,
\]
which means that $f_2$ is constant. Similarly, we get $f_1$ is constant. In this case, (\ref{Im1})  and
(\ref{Im2}) reduce to the following
\begin{equation}
\frac{\partial^2F_1^2}{\partial x^k\partial
y^l}y^k=\frac{\partial F_1^2}{\partial x^l},\label{Im3}
\end{equation}
\begin{equation}
\frac{\partial^2F_2^2}{\partial u^\alpha\partial
v^\beta}v^\alpha=2\frac{\partial F_2^2}{\partial
u^\beta}, \label{Im4}
\end{equation}
Hence $F_1$ and $F_2$ are locally dually flat.
\end{proof}

\smallskip

By Theorem \ref{THM1}, we conclude the following.

\begin{cor}
There is no  a locally dually flat proper DWP-Finsler manifold.
\end{cor}



\noindent
Esmaeil Peyghan\\
Department of Mathematics, Faculty  of Science\\
Arak University\\
Arak 38156-8-8349.  Iran\\
Email: epeyghan@gmail.com

\bigskip

\noindent
Akbar Tayebi\\
Department of Mathematics, Faculty  of Science\\
Qom University\\
Qom. Iran\\
Email: akbar.tayebi@gmail.com

\bigskip

\noindent
Behzad Najafi\\
Department of Mathematics, Faculty  of Science\\
Shahed University\\
Tehran. Iran\\
Email:\ najafi@shahed.ac.ir


\begin{thebibliography}{MaHo}
\bibitem{A} D. E. Allison, {\it Pseudoconvexity in Lorentzian doubly warped product}, Geom. Dedicata.  {\bf 39}(1991), 223-227.
\bibitem{am} S.-I. Amari, {\it Differential-Geometrical Methods in Statistics}, Springer Lecture Notes in Statistics,  Springer-Verlag, 1985.
\bibitem{amna} S.-I. Amari and H. Nagaoka, {\it Methods of Information Geometry}, AMS Translation of Math. Monographs,  Oxford University Press, 2000.
\bibitem{As} G. S. Asanov, {\it Finslerian extensions of Schwarzschild metric}, Fortschr. Phys. {\bf 40} (1992),
667-693.
\bibitem{As1} G. S. Asanov,  {\it Finslerian metric functions over the product $R\times M$ and their potential
applications}, Rep. Math. Phys. {\bf 41} (1998), 117-132.
 \bibitem{BM2} S.  B\'{a}cs\'{o} and M. Matsumoto, {\it On Finsler spaces of Douglas type, A generalization of notion of Berwald space}, Publ. Math. Debrecen. \textbf{51}(1997), 385-406.
\bibitem{Ba} D. Bao, {\it On two curvature-driven problems in Riemann-Finsler geometry}, Advanced
Studies in Pure Mathematics XX, 2007.
\bibitem {BCS} D. Bao, S.S. Chern and Z. Shen, {\it An Introduction to Riemannian-Finsler Geometry}, Spinger-Verlag,  2000.
\bibitem{BEP} J. K. Beem, P. Ehrlich and T. G. Powell, {\it Warped product manifolds in relativity}, in: Selected Studies: A Volume Dedicated to the Memory of Albert Einstein, (North-Holland, Armsterdarm), (1982),  41-56.
\bibitem{BP} J. B. Beem and T. G. Powell, {\it Geodesic completeness and maximality in Lorentzian warped
products}, Tensor (N. S.),  {\bf 39}(1982), 31-36.
\bibitem{ChSh} X. Chen and Z. Shen, {\it  On Douglas Metrics},  Publ. Math. Debrecen.  \textbf{66} (2005), 503-512.
\bibitem{G} A. Gebarowski, {\it On conformally recurrent doubly warped products}, Tensor (N. S.), {\bf 57} (1996), 192-196.
\bibitem{Ich}  Y. Ichijy\={o},  {\it Finsler spaces modeled on a Minkowski space}, J. Math. Kyoto Univ. {\bf 16} (1976), 639-652.
\bibitem{Koz} L. Kozma, I. R. Peter and C. Varga, {\it Warped product of Finsler-manifolds}, Ann. Univ. Sci. Pudapest 44,  {\bf 157} (2001).
\bibitem{Marian} M. I. Munteanu, {\it Doubly Warped Product CR-Submanifolds in Locally Conformal Kaehler Manifolds}, Monatshefte fur Mathematik, {\bf 150} (2007), 333-342.
\bibitem{Marian1} M. I. Munteanu, {\it A Note on Doubly Warped Product Contact CR-Submanifolds in trans-Sasakian Manifolds}, Acta Mathematica Hungarica, {\bf 116} (2007), 121-126.

\bibitem{NST1} B. Najafi, Z. Shen and A. Tayebi,  {\it On a projective class of Finsler metrics}, Publ. Math. Debrecen. {\bf 70}(2007), 211-219.
\bibitem{NST2} B. Najafi, Z. Shen and A. Tayebi, {\it Finsler metrics of scalar flag curvature with special non-Riemannian curvature properties}, Geom. Dedicata. {\bf 131} (2008), 87-97.
\bibitem{ShDiff} Z. Shen, {\it Differential Geometry of Spray and Finsler Spaces}, Kluwer Academic Publishers,  2001.
\bibitem{shenInf} Z. Shen, {\it Riemann-Finsler geometry with applications to information geometry}, Chin. Ann. Math. {\bf 27} (2006), 73-94.
\bibitem{TP} A. Tayebi and E. Peyghan, {\it On Ricci tensors of Randers metrics}, J. Geom. Phys. {\bf 60}(2010), 1665-1670.
\bibitem{U} B. \"{U}nal, {\it Doubly Warped Products}, Diff. Geom. Appl. {\bf 15} (3)(2001), 253-263.
\end{thebibliography}
\end{document}